# A note on the asymptotic distribution of the minimum density power divergence estimator

**Sergio F. Juárez[1] and William R. Schucany[2]**

*Veracruzana University and Southern Methodist University*

**Abstract:** We establish consistency and asymptotic normality of the minimum density power divergence estimator under regularity conditions different from those originally provided by Basu et al.

## 1. Introduction

Basu et al. [1] and [2] introduce the *minimum density power divergence estimator* (MDPDE) as a parametric estimator that balances infinitesimal robustness and asymptotic efficiency. The MDPDE depends on a tuning constant $\alpha \geq 0$ that controls this trade-off. For $\alpha = 0$ the MDPDE becomes the maximum likelihood estimator, which under certain regularity conditions is asymptotically efficient, see chapter 6 of Lehmann and Casella [5]. In general, as $\alpha$ increases, the robustness (bounded influence function) of the MDPDE increases while its efficiency decreases. Basu et al. [1] provide sufficient regularity conditions for the consistency and asymptotic normality of the MDPDE. Unfortunately, these conditions are not general enough to establish the asymptotic behavior of the MDPDE in more general settings. Our objective in this article is to fill this gap. We do this by introducing new conditions for the analysis of the asymptotic behavior of the MDPDE.

The rest of this note is organized as follows. In Section 2 we briefly describe the MDPDE. In Section 3 we present our main results for proving consistency and asymptotic normality of the MDPDE. Finally, in Section 4 we make some concluding comments.

## 2. The MDPDE

Let $G$ be a distribution with support $\mathcal{X}$ and density $g$. Consider a parametric family of densities $\{f(x;\theta) : \theta \in \Theta\}$ with $x \in \mathcal{X}$ and $\Theta \subseteq \mathbb{R}^p$, $p \geq 1$. We assume this family is identifiable in the sense that if $f(x;\theta_1) = f(x;\theta_2)$ a.e. in $x$ then $\theta_1 = \theta_2$. The *density power divergence* (DPD) between an $f$ in the family and $g$ is defined as

$$d_\alpha(g, f) = \int_\mathcal{X} \left\{ f^{1+\alpha}(x;\theta) - \left(1 + \frac{1}{\alpha}\right) g(x) f^\alpha(x;\theta) + \frac{1}{\alpha} g^{1+\alpha}(x) \right\} dx$$

[1]Facultad de Estadística e Informática, Universidad Veracruzana, Av. Xalapa esq. Av. Avila Camacho, CP 91020 Xalapa, Ver., Mexico, e-mail: sejuarez@uv.mx

[2]Department of Statistical Science, Southern Methodist University, PO Box 750332 Dallas, TX 75275-0332, USA, e-mail: schucany@smu.edu

for positive $\alpha$, and for $\alpha = 0$ as

$$d_0(g, f) = \lim_{\alpha \to 0} d_\alpha(g, f) = \int_{\mathcal{X}} g(x) \log[g(x)/f(x;\theta)] \mathrm{d}x.$$

Note that when $\alpha = 1$, the DPD becomes

$$d_1(g, f) = \int_{\mathcal{X}} [g(x) - f(x;\theta)]^2 \mathrm{d}x.$$

Thus when $\alpha = 0$ the DPD is the Kullback–Leibler divergence, for $\alpha = 1$ it is the $L^2$ metric, and for $0 < \alpha < 1$ it is a smooth bridge between these two quantities. For $\alpha > 0$ fixed, we make the fundamental assumption that there exists a unique point $\theta_0 \in \Theta$ corresponding to the density $f$ *closest* to $g$ according to the DPD. The point $\theta_0$ is defined as the target parameter. Let $X_1, \ldots, X_n$ be a random sample from $G$. The minimum density power estimator (MDPDE) of $\theta_0$ is the point that minimizes the DPD between the probability mass function $\hat{g}_n$ associated with the empirical distribution of the sample and $f$. Replacing $g$ by $\hat{g}_n$ in the definition of the DPD, $d_\alpha(g, f)$, and eliminating terms that do not involve $\theta$, the MDPDE $\hat{\theta}_{\alpha,n}$ is the value that minimizes

$$\int_{\mathcal{X}} f^{1+\alpha}(x;\theta) \mathrm{d}x - \left(1 + \frac{1}{\alpha}\right) \frac{1}{n} \sum_{i=1}^n f^\alpha(X_i;\theta)$$

over $\Theta$. In this parametric framework the density $f(\cdot;\theta_0)$ can be interpreted as the *projection* of the true density $g$ on the parametric family. If, on the other hand, $g$ is a member of the family then $g = f(\cdot;\theta_0)$.

Consider the score function and the information matrix of $f(x;\theta)$, $S(x;\theta)$ and $i(x;\theta)$, respectively. Define the $p \times p$ matrices $K_\alpha(\theta)$ and $J_\alpha(\theta)$ by

(2.1) $$K_\alpha(\theta) = \int_{\mathcal{X}} S(x;\theta) S^t(x;\theta) f^{2\alpha}(x;\theta) g(x) \mathrm{d}x - U_\alpha(\theta) U_\alpha^t(\theta),$$

where

$$U_\alpha(\theta) = \int_{\mathcal{X}} S(x;\theta) f^\alpha(x;\theta) g(x) \mathrm{d}x$$

and

(2.2) $$\begin{aligned} J_\alpha(\theta) &= \int_{\mathcal{X}} S(x;\theta) S^t(x;\theta) f^{1+\alpha}(x;\theta) \mathrm{d}x \\ &\quad + \int_{\mathcal{X}} \{i(x;\theta) - \alpha S(x;\theta) S^t(x;\theta)\} \times [g(x) - f(x;\theta)] f^\alpha(x;\theta) \mathrm{d}x. \end{aligned}$$

Basu et al. [1] show that, under certain regularity conditions, there exists a sequence $\hat{\theta}_{\alpha,n}$ of MDPDEs that is consistent for $\theta_0$ and the asymptotic distribution of $\sqrt{n}(\hat{\theta}_{\alpha,n} - \theta_0)$ is multivariate normal with mean vector zero and variance-covariance matrix $J_\alpha(\theta_0)^{-1} K_\alpha(\theta_0) J_\alpha(\theta_0)^{-1}$. The next section shows this result under assumptions different from those of Basu et al. [1].

## 3. Asymptotic Behavior of the MDPDE

Fix $\alpha > 0$ and define the function $m : \mathcal{X} \times \Theta \to \mathbb{R}$ as

(3.1) $$m(x,\theta) = \left(1 + \frac{1}{\alpha}\right) f^\alpha(x;\theta) - \int_{\mathcal{X}} f^{1+\alpha}(x;\theta) \mathrm{d}x$$

for all $\theta \in \Theta$. Then the MDPDE is an M-estimator with criterion function given by (3.1) and it is obtained by *maximizing*

$$m_n(\theta) = \frac{1}{n} \sum_{i=1}^{n} m(X_i, \theta)$$

over the parameter space $\Theta$. Let $\Theta_G \subseteq \Theta$ be the set where

(3.2) $$\int_{\mathcal{X}} |m(x,\theta)| g(x) \mathrm{d}x < \infty.$$

Clearly $\theta_0 \in \Theta_G$, but we assume $\Theta_G$ has more points besides $\theta_0$. For $\theta \in \Theta_G$ consider the expected value of $m(X;\theta)$ in (3.1) under the true distribution $G$

(3.3) $$M(\theta) = \left(1 + \frac{1}{\alpha}\right) \int_{\mathcal{X}} f^\alpha(x;\theta) g(x) \mathrm{d}x - \int_{\mathcal{X}} f^{1+\alpha}(x;\theta) \mathrm{d}x,$$

and define $M(\theta) = -\infty$ for $\theta \in \Theta \setminus \Theta_G$. Then the target parameter $\theta_0$ is such that $-\infty < M(\theta_0) = \sup_{\theta \in \Theta} M(\theta) < \infty$. Furthermore, we assume that $\Theta$ may be endowed with a metric $d$. Heretofore it is assumed that $(\Theta, d)$ is compact. The next theorem establishes consistency of the MDPDE.

**Theorem 1.** *Suppose the following conditions hold.*

1. *The target parameter $\theta_0 = \arg\max M_{\theta \in \Theta}(\theta)$ exists and is unique.*
2. *For $\theta \in \Theta_G$, $\theta \mapsto m(x,\theta)$ is upper semicontinuous a.e. in $x$.*
3. *For all sufficiently small balls $B \subset \Theta$, $x \mapsto \sup_{\theta \in B} m(x,\theta)$ is measurable and satisfies*

$$\int_{\mathcal{X}} \sup_{\theta \in B} m(x,\theta) g(x) dx < \infty.$$

*Then any sequence of MDPDEs $\hat{\theta}_{\alpha,n}$ that satisfies $m_n(\hat{\theta}_{\alpha,n}) \geq m_n(\theta_0) - o_p(1)$, is such that for any $\epsilon > 0$ and every compact set $K \subset \Theta$,*

$$P(d(\hat{\theta}_{\alpha,n}, \theta_0) \geq \epsilon, \hat{\theta}_{\alpha,n} \in K) \to 0.$$

*Proof.* This is Theorem 5.14 of van der Vaart [6] page 48. □

The first condition is our assumption of existence of $\theta_0$. It states that $\theta_0$ is an element of the parameter space and it is unique (identifiable). Without this assumption there is no minimum density power estimation to do. Compactness of $K$ is needed for $\{\theta \in K : d(\theta, \theta_0) \geq \epsilon\}$ to be compact; this is a technical requirement to prove the theorem. If $\Theta$ is not compact, one possibility is to compactify it. The third condition would follow if $f(x;\theta)$ is upper semicontinuous (trivially if it is continuous) in $\theta$ a.e. in $x$. Finally, the fourth condition is warranted by (3.2) in the interior of $\Theta_G$. Thus we can claim the following result.

**Theorem 2.** *If condition 1 in Theorem 1 holds, and if $f^\alpha(x;\theta)$ is upper semicontinuous (continuous) in $\theta$ in the interior of $\Theta_G$ and for a.e. in $x$, then any sequence $\hat{\theta}_{\alpha,n}$ of MDPDEs such that $m_n(\hat{\theta}_{\alpha,n}) \geq m_n(\theta_0) - o_p(1)$, satisfies $d(\hat{\theta}_{\alpha,n}, \theta_0) \xrightarrow{p} 0$.*

The asymptotic normality of the MDPDE hinges on smoothness conditions that are not required for consistency. These conditions are provided in the two following results.

**Lemma 3.** $M(\theta)$ *as given by* (3.3) *is twice continuous differentiable in a neighborhood $B$ of $\theta_0$ with second derivative (Hessian matrix)* $\mathsf{H}_\theta M(\theta) = -(1+\alpha)J_\alpha(\theta)$, *if:*

1. *The integral $\int_\mathcal{X} f^{1+\alpha}(x;\theta)dx$ is twice continuously differentiable with respect to $\theta$ in $B$, and the derivative can be taken under the integral sign.*
2. *The order of integration with respect to $x$ and differentiation with respect to $\theta$ can be interchanged in $M(\theta)$, for $\theta \in B$.*

*Proof.* Consider the (transpose) score function $S^t(x;\theta) = \mathsf{D}_\theta \log f(x;\theta)$ and the information matrix $i(x;\theta) = -\mathsf{H}_\theta \log f(x;\theta) = -\mathsf{D}_\theta S(x;\theta)$. Also note that $[\mathsf{D}_\theta f(x;\theta)]f^{\alpha-1}(x;\theta) = S^t(x;\theta)f^\alpha(x;\theta)$. Use the previous expressions and condition 1 to obtain the first derivative of $\theta \mapsto m(x;\theta)$

$$(3.4) \quad \mathsf{D}_\theta m(x,\theta) = (1+\alpha)S^t(x;\theta)f^\alpha(x;\theta) - (1+\alpha)\int_\mathcal{X} S^t(x;\theta)f^{1+\alpha}(x;\theta)\mathrm{d}x.$$

Proceeding in a similar way, the second derivative of $\theta \mapsto m(x;\theta)$ is

$$(3.5) \quad \begin{aligned}\mathsf{H}_\theta m(x,\theta) &= (1+\alpha)\{-i(x;\theta) + \alpha S(x;\theta)S^t(x;\theta)\}f^\alpha(x;\theta) - (1+\alpha) \\ &\quad \times \left\{\int_\mathcal{X} -i(x;\theta)f^{1+\alpha}(x;\theta) + (1+\alpha)S(x;\theta)S^t(x;\theta)f^{1+\alpha}(x;\theta)\mathrm{d}x\right\}.\end{aligned}$$

Then using condition 2 we can compute the second derivative of $M(\theta)$ under the integral sign and, after some algebra, obtain

$$\mathsf{H}_\theta M(\theta) = \int_\mathcal{X} \{\mathsf{H}_\theta m(x,\theta)\}g(x)\mathrm{d}x = -(1+\alpha)J_\alpha(\theta). \qquad \square$$

The second result is an elementary fact about differentiable mappings.

**Proposition 4.** *Suppose the function $\theta \mapsto m(x,\theta)$ is differentiable at $\theta_0$ for $x$ a.e. with derivative $\mathsf{D}_\theta m(x,\theta)$. Suppose there exists an open ball $B \in \Theta$ and a constant $M < \infty$ such that $\|\mathsf{D}_\theta m(x,\theta)\| \leq M$ for all $\theta \in B$, where $\|\cdot\|$ denotes the usual Euclidean norm. Then for every $\theta_1$ and $\theta_2$ in $B$ and a.e. in $x$, there exist a constant that may depend on $x$, $\phi(x)$, such that*

$$(3.6) \quad |m(x,\theta_1) - m(x,\theta_2)| \leq \phi(x)\|\theta_1 - \theta_2\|,$$

*and*

$$\int_\mathcal{X} \phi^2(x)g(x)dx < \infty.$$

We can now establish the asymptotic normality of the MDPDE.

**Theorem 5.** *Let the target parameter $\theta_0$ be an interior point of $\Theta$, and suppose the conditions of Lemma 3 and Proposition 4 hold. Then, any sequence of MDPDEs $\hat\theta_{\alpha,n}$ that is consistent for $\theta_0$ is such that*

$$\sqrt{n}(\hat\theta_{n,\alpha} - \theta_0) \rightsquigarrow N_p(0, J_\alpha^{-1}(\theta_0)K_\alpha(\theta_0)J_\alpha^{-1}(\theta_0)),$$

*where $K_\alpha$ and $J_\alpha$ are given in* (2.1) *and* (2.2), *respectively.*

*Proof.* From Lemma 3, $M(\theta)$ admits the following expansion at $\theta_0$

$$M(\theta) = M(\theta_0) + \frac{1}{2}(\theta - \theta_0)^t V_\alpha(\theta_0)(\theta - \theta_0) + o(\|\theta - \theta_0\|^2),$$

where $V_\alpha(\theta) = \mathsf{H}_\theta M(\theta)$. Proposition 4 implies the Lipschitz condition (3.6). Then the conclusion follows from Theorem 5.23 of van der Vaart [6]. $\square$

So far we have not given explicit conditions for the existence of the matrices $J_\alpha$ and $K_\alpha$ as defined by (2.2) and (2.1), respectively. In order to complete the asymptotic analysis of the MDPDE we now do that. Condition 2 in Lemma 3 implicitly assumes the existence of $J_\alpha$. This can be justified by observing that the condition that allows interchanging the order integration and differentiation in $M(\theta)$ is equivalent to the existence of $J_\alpha$. For $J_\alpha$ to exist we need $i_{jk}(x;\theta)$, the $jk$-element of the information matrix $i(x;\theta)$, to be such that

$$\int_\mathcal{X} i_{jk}(x;\theta) f^{1+\alpha}(x;\theta) \mathrm{d}x < \infty, \text{ and } \int_\mathcal{X} i_{jk}(x;\theta) f^\alpha(x;\theta) g(x) \mathrm{d}x < \infty.$$

Regarding $K_\alpha$, let $S_j(x;\theta)$ be the $j$th component of the score $S(x;\theta)$. If

(3.7) $$S_j^2(x;\theta) f^{2\alpha}(x;\theta) < C_j, \quad j = 1, \ldots, p,$$

then

$$\int_\mathcal{X} S_j^2(x;\theta) f^{2\alpha}(x;\theta) g(x) \mathrm{d}x < \infty.$$

Thus, $K_\alpha$ exists. Furthermore, by (3.4) we see that the $j$th component of $\mathsf{D}_\theta m(x,\theta)$ is

$$A_j = (1+\alpha) S_j(x;\theta) f^\alpha(x;\theta) - (1+\alpha) \int_\mathcal{X} S_j(x;\theta) f^{1+\alpha}(x;\theta) \mathrm{d}x.$$

Then $A_j^2$ would be bounded by a constant $M_j < \infty$ if all the components $S_j(x;\theta)$ of the score vector $S(x;\theta)$ satisfy (3.7). This is true because in this case $S_j(x;\theta) f^\alpha(x;\theta)$ would be bounded by a constant too, and then

$$\int_\mathcal{X} S(x;\theta) f^{1+\alpha}(x;\theta) \mathrm{d}x < \infty.$$

Hence

$$\| \mathsf{D}_\theta m(x,\theta) \|^2 = \sum_{i=1}^p A_j^2 \leq \sum_{i=1}^p M_j < \infty.$$

Therefore, if (3.7) holds, then the Lipschitz condition in (3.6) follows.

From the previous analysis, we see that the conditions on $M(\theta)$ and $m(x,\theta)$ given in Theorem 5 can be established in terms of the density $f(x;\theta)$, its score vector $S(x;\theta)$, and its information matrix $i(x;\theta)$ as indicated in the next theorem.

**Theorem 6.** *The MDPDE is asymptotically normal, as in Theorem 5, if the following conditions hold in a neighborhood $B \subseteq \Theta_G$ of $\theta_0$:*

1. *Condition 1 of Lemma 3*
2. *For each $j = 1, \ldots, p$, $S_j^2(x;\theta) f^{2\alpha}(x;\theta) < C_j$ a.e. in $x$.*
3. *Suppose there are functions $\phi_{jk}$ such that $|i_{jk}(x;\theta) f^\alpha(x;\theta)| \leq \phi_{jk}(x)$ for $j, k = 1, \ldots, p$, and*

$$\int_\mathcal{X} \phi_{jk}(x) f(x;\theta) dx < \infty, \text{ and } \int_\mathcal{X} \phi_{jk}(x) g(x) dx < \infty.$$

## 4. Concluding Remarks

We have obtained consistency of the MDPDE under rather general conditions on the criterion function $m$. Namely, integrability of $x \mapsto m(x,\theta)$, and upper semicontinuity of $\theta \mapsto m(x,\theta)$. However these are not necessary conditions; concavity or asymptotic concavity of $m_n(\theta)$, would also give consistency of the MDPDE without requiring compactness of $\Theta$, see Giurcanu and Trindade [3]. To decide which set of conditions are easier to verify seems to be more conveniently handled on a case by case basis.

**Acknowledgements**

The authors thank professor Javier Rojo for the invitation to present this work at the Second Symposium in Honor of Erich Lehmann held at Rice University. They are also indebted to the editor for his comments and suggestions which led to a substantial improvement of the article. Finally, the first author is deeply grateful to Professor Rojo for his proverbial patience during the preparation of this article.